\magnification=\magstephalf
\input diagrams

\font\tengoth=eufm10  \font\fivegoth=eufm5
\font\sevengoth=eufm7
\newfam\gothfam  \scriptscriptfont\gothfam=\fivegoth
\textfont\gothfam=\tengoth \scriptfont\gothfam=\sevengoth
\def\goth{\fam\gothfam\tengoth}
%
\newfam\bboldfam
\font\tenbbold=msbm10
\font\sevenbbold=msbm7
\font\fivebbold=msbm5
\textfont\bboldfam=\tenbbold
\scriptfont\bboldfam=\sevenbbold
\scriptscriptfont\bboldfam=\fivebbold
\def\bbold{\fam\bboldfam\tenbbold}

\newcount\cols
{\catcode`,=\active\catcode`|=\active
\gdef\Young(#1){\hbox{$\vcenter
{\mathcode`,="8000\mathcode`|="8000
\def,{\global\advance\cols by 1 &}%
\def|{\cr
      \multispan{\the\cols}\hrulefill\cr
       &\global\cols=2 }%
  \offinterlineskip\everycr{}\tabskip=0pt
  \dimen0=\ht\strutbox \advance\dimen0 by \dp\strutbox
    \halign
    {\vrule height \ht\strutbox depth \dp\strutbox##
      &&\hbox to \dimen0{\hss$##$\hss}\vrule\cr
     \noalign{\hrule}&\global\cols=2 #1\crcr
     \multispan{\the\cols}\hrulefill\cr%
   }
}$}}
}
\def\sqr#1#2{{\vcenter{\vbox{\hrule height.#2pt
\hbox{\vrule width.#2pt height#1pt \kern#1pt
\vrule width.#2pt}
\hrule height.#2pt}}}}

\def \Box {\hfill\hbox{}\nobreak \vrule width 1.6mm height 1.6mm}

\def \fix#1 {{\hfill\break \bf ((*** #1 ***))\hfill\break}}

\def\g{{\goth g}}
\def\gl{{\goth gl}}

\def\Tor{{\rm Tor}}
\def\Ext{{\rm Ext}}
\def\S{{\cal S}}
\def\F{{\cal F}}
\def\FF{{\bf F}}
\def\th{{\rm th}}

\def\coker{\mathop{\rm coker}\nolimits}
\def\ann{\mathop{\rm ann}\nolimits}

\def\NN{{\bbold N}}
\def\ZZ{{\bbold Z}}

\centerline {\bf A Fitting Lemma for $\ZZ/2$-graded modules}
\medskip
\centerline {by}
\medskip
\centerline{\bf  David Eisenbud and Jerzy Weyman%
\footnote*{\rm The second named author is grateful to the
Mathematical Sciences Research Institute for support in the period this
work was completed. Both authors are grateful
for the partial support of the National Science Foundation.
{\bf Version of 2/18/02}}}
\medskip

\noindent {\bf Abstract.}
{\narrower { 
We study the annihilator of the cokernel of a map of free
$\ZZ/2$-graded modules over a $\ZZ/2$-graded skew-commutative algebra
in characteristic 0 and define analogues of its Fitting ideals.  We
show that in the ``generic'' case the annihilator is given by a
Fitting ideal, and explain relations between the Fitting ideal and the
annihilator that hold in general.  Our results generalize the
classical Fitting Lemma, and extend the key result of Green [1999].
They depend on the Berele-Regev theory of representations
of general linear Lie super-algebras.
}
\medskip
}
\beginsection Introduction.

The classical Fitting Lemma (Fitting [1936]) gives information
about the annihilator of a module over a commutative ring in terms
of a presentation of the madule by generators and relations.
More precisely, let
$$
\phi: R^m \to R^d
$$
be a map of finitely generated free modules over a commutative ring
$R$, and for any integer $t\geq 0$ 
let $I_t(\phi)$ denote the ideal in $R$ generated by the 
$t\times t$ minors of $\phi$.
Fitting's result says that the module $\coker \phi$
is annihilated by $I_d(\phi)$, and that if $\phi$ is the 
generic map---represented by a matrix whose entries are distinct indeterminates---then the
annihilator is equal to $I_d(\phi)$. Thus $I_d(\phi)$ is the best approximation to the
annihilator that is compatible with base change. Morever, 
$I_d(\phi)$ is not too bad an approximation to $\ann \coker \phi$
in the sense that 
$I_d(\phi)\supset (\ann \coker \phi)^d$, or more precisely
$(\ann\coker \phi)I_t(\phi)\subset I_{t+1}(\phi)$ for all $0\leq t<d$.

In this paper we will
prove a corresponding result in the case of $\ZZ/2$-graded modules over a
skew-commutative $\ZZ/2$-graded algebra containing a field $K$
of characteristic 0.
Let $R$ be a 
$\ZZ/2$-graded skew-commutative $K$-algebra: that is, $R=R_0\oplus
R_1$ as vector spaces, $R_0$ is a commutative central subalgebra,
$R_iR_j\subset R_{i+j\ \rm (mod\ 2)}$, and every element of
$R_1$ squares to 0. Any homogeneous map $\phi$ of $\ZZ/2$-graded free
$R$-modules may be written in the form
$$
\phi:\quad R^m\oplus R^n(1) \rTo^{\pmatrix{X&A\cr B&Y}} R^d\oplus R^e(1),
$$
where $X, Y$ are matrices of even elements of $R$, and $A, B$ are
matrices of odd elements. We will define an ideal $I_{\Lambda(d,e)}$ and
show that it is contained in the
annihilator of
the cokernel of $\phi$, with equality 
in the {\it generic case\/} where the entries of the matrices
$X,Y,A,B$ are
indeterminates (that is,
$R$ is a polynomial ring on the entries of $X$ and $Y$
tensored with an exterior algebra on the entries of $A$ and $B$.) We 
give examples to show that the annihilator can be quite different in positive
characteristic.
\goodbreak

Now let $K$ be a field, and let 
$U=U_0\oplus U_1$ and $V=V_0\oplus V_1$ be $\ZZ/2$-graded vector
spaces of dimensions
$(d,e)$ and $(m,n)$ respectively.
We consider the generic ring
$$\S=
\S(V\otimes U) := S(V_0\otimes U_0 )\otimes S(V_1\otimes U_1
)\otimes\wedge (V_0\otimes U_1
)\otimes\wedge(V_1\otimes U_0 ),
$$
where $S$ denotes the symmetric algebra and $\wedge$ denotes
the exterior algebra, and the {\it generic\/}, or
{\it tautological \/} map
$$
\Phi: \S\otimes V\rTo \S\otimes U^*.
$$
This map $\Phi$ is defined
by the condition that
$\Phi|_V=1\otimes \eta: V\to V\otimes_K U\otimes_K U^*\subset R\otimes_K U^*$,
where $\eta: K\to U\otimes_K U^*$ is the dual of the contraction
$U^*\otimes_K U\to K$.
We will make use of this notation throughout the paper. 

We will
compute the annihilator of the cokernel of $\Phi$. Of course if
we specialize $\Phi$ to any map of free modules $\phi$ over a
$\ZZ/2$-graded ring,
preserving the grading, then we can derive elements in the
annihilator of the cokernel of $\phi$ by specializing the
annihilator of the cokernel of $\Phi$.

In the classical case, where $V$ and $U$ have only even parts
($e=n=0$) the annihilator is an invariant ideal for the action of the
product of general linear groups ${\rm GL}(V)\times {\rm GL}(U)$.
Such invariant ideals have been studied by DeConcini, Eisenbud, and
Procesi in [1980] and have a very simple arithmetic.  In the general
case Berele and Regev [1987] have developed a highly parallel theory,
using the $\ZZ/2$-graded Lie algebra $\g={\gl}(V)\times {\gl}(U)$ in
place of ${\rm GL}(V)\times {\rm GL}(U)$.  They show that the generic
ring $\S$ is a semisimple representation of $\g$ (even though not all
the representations of $\g$ are semisimple)  and that the irreducible
summands of $\S$ of total degree $t$ are parametrized by certain
partitions of the integer $t$, just as in the commutative case.
The Berele-Regev theory is described in detail below, in Section 1 of this
paper.  

If $\Lambda$ is a partition, we write $I_\Lambda$ for the
ideal of $\S$ generated by the irreducible representation
corresponding to $\Lambda$. If $\phi$ is a matrix representing any map
of $\ZZ/2$-graded free modules over a $\ZZ/2$-graded skew-commutative
$K$ algebra $R$, then there is a unique ring homomorphism $\alpha:
\S\to R$ such that $\phi = \alpha(\Phi)$, and we write
$I_\Lambda(\phi):=\alpha(I_\Lambda(\Phi))R$ for the ideal generated by
the image of $I_\Lambda=I_\Lambda(\Phi)$.

If $e=n=0$ the classical Fitting
Lemma shows that the annihilator of the module with generic presentation
matrix as above is $I_d(\Phi)$. In representation-theoretic terms
this is the ideal 
generated by the representation $\wedge^dV\otimes \wedge^dU$, 
the irreducible representation associated to the partition with
one term $(d)$. In our notation, $I_d(\Phi)= I_{(d)} = I_{(d)}(\Phi)$). Here is the
generalization which is the main result of this paper:

\proclaim Theorem 1. Suppose that $K$ is a field of characteristic 0,
and let  
$$
\phi:\quad R^m\oplus R^n(1) \rTo^{\pmatrix{X&A\cr B&Y}} R^d\oplus R^e(1),
$$
be a $\ZZ/2$-graded map of free modules over a $\ZZ/2$-graded skew-commutative
$K$-algebra $R$.
\hfill\break\vskip .1cm
\noindent{\bf a) } When $R=\S$ and $\phi=\Phi$, the generic map
defined above,
the annihilator of the cokernel of $\Phi$ is
$
I_{\Lambda(d,e)}(\Phi),
$
where
$\Lambda(d,e)$ is the partition $(d+1,d+1,\dots,d+1, d)$ of $(d+1)(e+1)-1$
into $e+1$ parts. In general we have
$I_{\Lambda(d,e)}(\phi)\subset \ann\ \coker (\phi)$.
\hfill\break\vskip .1cm
\noindent{\bf b)} 
If $x_1 ,\ldots ,x_e \in
\ann\ \coker(\phi)$, then
$x_1\ldots x_e\in I_{\Lambda (0,e)}(\phi)$.
Moreover, if  $0\le s\le d-1$, and 
$x_1 ,\ldots ,x_{e+1}\in \ann\ \coker(\phi)$, 
then 
$
x_1\ldots x_{e+1}I_{\Lambda (s,e)}(\phi )\subset I_{\Lambda (s+1,e)}(\phi ).
$

\noindent The proof is given in sections 2 and 3 below.

In the classical case ($e=n=0$) we can also describe the
annihilator of $\coker \Phi$ by saying that it is nonzero only if
$m\geq d$, and then it is generated, as a
${\gl}(V)\times {\gl}(U)$-ideal, by an
$m\times m$ minor of $\Phi$. To simplify the general statement
we note that a shift of degree
by 1 does not change the annihilator of the cokernel of $\Phi$,
but has the effect of interchanging $m$ with $n$ and $d$ with $e$.

\proclaim Corollary 2. With notation as above, 
the annihilator of the cokernel of $\Phi$ is
nonzero only if\hfill\break
{a)} $m> d$ (or symmetrically $n> e$) or\hfill\break
{b)} $m=d$ and $n=e$.\hfill\break
\noindent In each of these cases the annihilator is generated as
a $\g$-ideal by one element $Z$ of degree $de+d+e$ defined as follows:
\hfill\break
\indent In case $a)$ when $m>d$
$$Z=Z_1 \cdot X(1,\ldots ,d|1,\ldots ,d )$$
where $X(1,\ldots ,d|1,\ldots
,d)$ is the
$d\times d$ minor of $X$
corresponding to the first $d$ columns and
$Z_1 =\prod_{j\leq e,k\leq d+1}b_{j,k}$ is the product
of all the elements in the first $d+1$ columns of $B$
(and symmetrically if $n>e$);\hfill\break
%
%
\indent In case $b)$
$$Z= W_1\cdots W_e\cdot det(X)$$
where $W_s$ is the $(d+1)\times (d+1)$ minor of
$\Phi$ containing $X$ and the entry $y_{s,s}$, that is,
$$W_s = det(X)y_{s,s}+\sum_{1\le i,\le d}\pm det(X({\hat i},{\hat
j}))a_{i,s}b_{s,j}.$$
\Box

\noindent Corollary 2 follows from intermediate results in the proof of Theorem
1a. We next give some examples of Theorem 1a and Corollary 2.
\medskip

\noindent{\bf Example 1.} Suppose that $d=n=0$, so that the presentation
matrix $B$ has only odd degree entries.
A central observation of Green [1999]
is that the ``exterior minors'' of $\Phi$ are in the annihilator of
$\coker \Phi$.
The element $Z$ of Corollary 2 is the product of the elements
in the first column of $\Phi$. Quite generally, it is not hard to see that the product
of all elements
in a $K$-linear combination of the columns of $B$ is an exterior
minor in Green's sense. 
The representation corresponding to the
partition $(1,\dots,1)$ of $e$ is generated by
${m+e-1\choose e}$ such products, so 
 ${m+e-1\choose e}$ exterior
minors generate the annihilator in the generic case.

For
example, taking $m=2,\ e=2$ the annihilator of the
cokernel of the generic matrix
$$\pmatrix{
b_{1,1}&b_{1,2}\cr
b_{2,1}&b_{2,2}
},
$$
where the variables all have odd degree, is minimally generated by the three
exterior minors
$$
b_{1,1}b_{2,1},\quad b_{1,2}b_{2,2},\quad
(b_{1,1}+b_{1,2})(b_{2,1}+b_{2,2}).
$$
\medskip

\noindent{\bf Example 2.} Now suppose that
our generic matrix has size $2\times 2$
with the first row even and the second row odd ($m=2,\; n=0,\; d=e=1$):
$$
\pmatrix{
x_{1,1}&x_{1,2}\cr
b_{1,1}&b_{1,2}
}.
$$
In this case our result shows that the cokernel has annihilator equal
to the product
$$
(x_{1,1},\quad x_{1,2})(b_{1,1}b_{1,2},\quad x_{1,1}b_{1,2}-x_{1,2}b_{1,1}),
$$
which is minimally generated by 4 elements. The element $Z$ is
$x_{1,1}b_{1,1}b_{1,2}$.
\medskip

\noindent{\bf Example 3.} As a final $2\times 2$ example, consider the case
$m=n=d=e=1$ which for simplicity we write as
$$
\pmatrix{
x&a\cr
b&y}.
$$
Here the annihilator of the cokernel is again minimally generated
by 4 elements, namely
$$
axy,\quad bxy,\quad (xy-ab)x,\quad (xy+ab)y.
$$
The element $Z$ is $(xy+ab)x$.
In an Appendix we will explain the $\g$ action on these elements.

\medskip

\noindent{\bf Positive characteristics.}
Already with $m=n=d=e=1$ as in Example 3 the annihilator is different 
in characteristic 2: in characteristic zero the annihilator is generated by forms of
degree 3, but in characteristic 2 the
algebra $R$ is commutative, so
the determinant $xy-ab$ is in the annihilator as well.

The annihilator can differ in other characteristics as well. 
Macaulay2 computations show that the case
$d=1,\ e=p-1,\ m=2,\ n=0$ is exceptional in characteristic $p$
for $p=3, 5$  and 7. Perhaps the same holds
for all primes $p$.

The cokernel of the generic matrix over the integers
can also have $\ZZ$-torsion. For example
Macaulay2 computation shows that if $d=1,\ e=2,\ m=3,\ n=1$
then the cokernel of $\Phi_\ZZ$ has 2-torsion.
\bigskip

Our interest in extending the Fitting Lemma
was inspired by Mark Green's paper
[1999] where he shows that the exterior minors are
in the annihilator. Green's striking use of
his result to prove one of the
Eisenbud, Koh, Stillman conjectures on
linear syzygies turns on the fact that if $N$ is a module
over a polynomial ring $S=K[X_1,\dots,X_m]$ then $T:=\Tor_*^S(K,M)$
is a module over the ring $R=\Ext_S^*(K,K)$, which is an
exterior algebra. Green in effect
translated the hypothesis of the linear syzygy
conjecture into a statement about the degree 1 part of the $R$-free
presentation
matrix of the submodule of $T$ representing the linear part of the
resolution of
$N$, and then showed that the exterior minors generated a certain
power of the maximal ideal of the exterior algebra, which was sufficient
to prove the Conjecture.
Green's result only gives information on the annihilator in the
case where the elements of the presentation matrix are all odd.
Elements of even degree in an exterior algebra can behave (if the
number of variables is large) very much like variables in a polynomial
ring, at least as far as expressions of bounded degree are concerned.
Thus to extend Green's work it seemed natural to deal with the
case of $\ZZ/2$-graded algebras.

This work is part of a program to study modules and resolutions over
exterior algebras; see Eisenbud-Fl\o ystad-Schreyer [2001], 
and Eisenbud-Popescu-Yuzvinsky [2000] for further information.

We would never have undertaken the project reported in this paper
if we had not had the program Macaulay2 (www.math.uiuc.edu/Macaulay2) of
Grayson and Stillman as a tool; its ability to compute in skew commutative
algebras was invaluable in figuring out the pattern that the results should
have and in assuring us that we were on the right track.

\beginsection 1. Berele-Regev Theory

For the proof of Theorem 1 we will use the beautiful results
of Berele and Regev [1987] giving the structure of $R$ as
a module over $\g$. For the convenience of the reader
we give a brief sketch of what is needed.
{\it
We make use of the notation introduced above:
$U=U_0\oplus U_1$ and $V=V_0\oplus V_1$ are $\ZZ/2$ graded vector spaces
over the field  $K$ of characteristic 0 with $\dim U=(d,e)$ and
$\dim V=(m,n)$.}

The
$\ZZ/2$-graded Lie algebra
$\gl(V)$ is the vector space of $\ZZ/2$-graded
 endomorphisms of $V=V_0\oplus V_1$. Thus
$$\gl(V)= \gl(V)_0\oplus \gl(V)_1,$$
where $\gl(V)_0$ is the set of endomorphisms preserving the grading of $V$
and $\gl(V)_1$ is the set of
endomorphisms of $V$ shifting the grading by 1. Additively
$$\gl(V)_0 = End_K (V_0 )\oplus End_K (V_1 ),$$
$$ \gl(V)_1 =Hom_K (V_0 ,V_1 )\oplus Hom_K (V_1 ,V_0 )$$
The commutator of the pair of homogeneous elements $x,y\in \gl(V)$ is
defined
by the formula
$$[x,y]= xy-(-1)^{deg(x)deg(y)}yx.$$

By
a $\gl(V)$-module we mean a $\ZZ/2$-graded vector space $M=M_0\oplus M_1$
with a bilinear map of $\ZZ/2$-graded vector spaces
$\circ :\gl(V)\times M\rightarrow M$   satisfying the identity
$$
[x,y]\circ m = x\circ (y\circ m)-(-1)^{deg(x)deg(y)} y\circ (x\circ m))
$$
for homogeneous elements $x,y\in \gl(V), m\in M$.

In contrast to the classical theory, not every representation of the
$\ZZ/2$-graded Lie algebra $\gl(V)$ is
semisimple. For example its natural action on mixed tensors $V^{\otimes
k}\otimes V^*{}^{\otimes l}$ is in
general not completely reducible. However, its action on
$V^{\otimes t}$ decomposes just as in the ungraded case:

\proclaim Proposition 1.1. The action of $\gl(V)$ on $V^{\otimes t}$ is
completely reducible for each $t$.
More precisely, the analogue of Schur's double centralizer theorem holds
and the irreducible $\gl(V)$-modules
occurring in the decomposition of $V^{\otimes t}$  are in 1-1
correspondance with irreducible representations of
the symmetric group $\Sigma_t$ on $t$ letters. These irreducibles are the Schur
functors
$$\S_\lambda (V) = e(\lambda )V^{\otimes t}$$
where $e(\lambda )$ is a Young idempotent corresponding to a partition
$\lambda$ in the group ring of the
symmetric group $\Sigma_t$.\Box

This notation is consistent with the notation above in the sense that the $d$-th
homogeneous component of the ring $\S (V)$ is $\S_d (V)$ where $d$ represents the
partition $(d)$ with one part.

Here we use the symbol $\S_\lambda$ to denote the $\ZZ/2$-graded
version of the Schur functor $S_\lambda$; the latter acts on ungraded vector
spaces. Recall that the functor $\wedge^\lambda$ is by definition the same
as the functor $S_{\lambda'}$, where $\lambda'$ denotes the partition
which is conjugate to $\lambda$. (For example, the
conjugate partition to $(2)$ is $(1,1)$.)
We will extend this by writing
$\bigwedge^\lambda := \S_{\lambda'}$ for the $\ZZ/2$-graded vesion.
The partition $(d)$ with only one part will be denoted simply $d$, so
for example
$\S_2(V)=\bigwedge^{(1,1)}V=S_2(V_0)\oplus (V_0\otimes
V_1)\oplus
\wedge^2V_1$ and similarly
$\bigwedge^{2}V=\S_{(1,1)}V=\wedge^2V_0\oplus V_0\otimes V_1\oplus S_2(V_1)$.
In each case the decomposition is as representations of the subalgebra
$\gl(V_0)\times \gl(V_1)\subset \gl(V)$. Similar decompositions hold for
all $\S_d$ and $\bigwedge^{d}V$. (If we were not working in
characteristic zero we would use divided powers in place of
symmetric powers in the description
of $\bigwedge V$.)

Proposition 1.1 implies that the
parts of the representation theory of $\gl(V)\times \gl(U)$
that involve only tensor products of
$V$ and $U$ and their summands are parallel to the representation
theory in the case $V_1 =U_1 =0$, which
is the classical representation theory of product
of the two general linear Lie algebras ${\rm gl}(V_0)\times {\rm gl}(U_0)$.

The Proposition also implies that the decompositions into irreducible
representations of tensor products of the $\S_\lambda(V)$, as well as the
decompositions of their symmetric
and exterior powers,  correspond to the
decompositions in the even case: we just have
to replace the ordinary Schur functors $S, \wedge$ by their
$\ZZ/2$-graded analogues $\S, \bigwedge$.

The
formulas giving equivariant embeddings or
equivariant projections may also be derived from the corresponding formulas
in the even case by
applying the principle of signs: The formulas in the even case involve many
terms where the basis elements are
permuted in a prescribed way. The basis elements have degree $0$. To write
down a $\ZZ/2$-graded analogue of
such formula we simply allow the basis elements to have even or odd degree
and we adjust the signs of terms in
such way that changing the order of two
homogeneous elements $x$ and $y$ of
$V$  in the $\ZZ/2$-graded analogue of the formula
will cost the additional factor
$(-1)^{deg(x)deg(y)}$.

There is a $\ZZ/2$-graded analogue of the Cauchy decomposition,
which follows as just described from Proposition 1.1 together with
the corresponding result in the even case (proven in [MD], ch.1
and [DC-E-P]). Recall that $\g=\gl(V)\times\gl(U)$.

\proclaim Corollary 1.2. The $t^\th$ component
$\S_t (V\otimes U)$ of $\S(V\otimes U)$
decomposes as a $\g$-module as
$$
\S_t (V\otimes U)=\oplus_{\lambda ,|\lambda |=t} \S_\lambda(V)\otimes
\S_\lambda (U).\quad \Box
$$

Another application of the same principle shows that to
describe the annihilator of the cokernel of $\Phi$, and what
generates it, it suffices to describe which
representations $\S_\lambda V\otimes
\S_\lambda U$ it contains:

\proclaim Corollary 1.3. If $I\subset \S(V\otimes U)$ is a $\g$-invariant
ideal, then $I$ is a sum of subrepresentations $\S_\lambda V\otimes
\S_\lambda U$.
Moreover, the ideal generated by $\S_\lambda V\otimes \S_\lambda U$ contains
$\S_\mu V\otimes \S_\mu U$ if and only if
$\mu\supset\lambda$.\Box

Although it is not so simple to describe the vectors in $\S_t(V\otimes U)$
that lie in a given irreducible summand, we can, as
in the commutative case, define a filtration that has these
irreducible representations as successive factors.
We start by defining a map
$\rho_t : \bigwedge^t V\otimes\bigwedge^t
U\hookrightarrow \S_t (V\otimes U)$ as the composite
$$
\bigwedge^t V\otimes\bigwedge^tU\rTo
\otimes^tV\ \otimes\ \otimes^tU\rTo
\S_t(V\otimes U)
$$
where the first map is the tensor product of the two diagonal maps
(here we use the sign conventions for $\ZZ/2$-graded vector spaces)
and the second map simply pairs corresponding factors. Thus
$$
\rho_t (v_1\wedge \ldots\wedge v_t\otimes u_1\wedge\ldots\wedge u_t
)=\sum_{\sigma\in\Sigma_t} \pm (v_1\otimes
u_{\sigma(1)})\cdot\ldots\cdot (v_t\otimes u_{\sigma(t)})
$$
 where the sign $\pm$ is the sign of the
permutation $\sigma$ adjusted by
the rule that switching homogeneous elements
$x, y$ from either $V$ or $U$ means we multiply by $(-1)^{deg(x)deg(y)}$.
For example, if $V$ and $U$
were both even, the image of this map would be
the span of the $t\times t$ minors of the generic matrix;
when $V$ is even and $U$ is odd, the image is the span of the
space of ``exterior minors" as in Green [1999].

For any partition $\lambda=(\lambda_1,\dots,\lambda_s)$
we define  $\F_\lambda$ to be
the image of the composite map
$$
m\circ(\rho_{\lambda_1}\otimes\ldots\otimes\rho_{\lambda_s}):\bigwedge^{\lambda_
1}
V\otimes\bigwedge^{\lambda_1}U
\otimes\ldots\otimes
\bigwedge^{\lambda_s}V\otimes\bigwedge^{\lambda_s}U\rightarrow \S_{|\lambda
| }(V\otimes U)
$$
where $m$ denotes the multiplication map in $\S(V\otimes U)$.

As in the even case we order partitions of $t$ by saying $\lambda < \mu$ if
and only if $\lambda^\prime_i
>\mu^\prime_i$ for the smallest number such that $\lambda^\prime_i
\ne
\mu^\prime_i$.
Finally, we define the subspaces
$$
\F_{<\lambda}=\sum_{|\mu|=|\lambda|,\ \mu<\lambda}\F_\mu
\qquad \subset \qquad
\F_{\le \lambda}=\sum_{|\mu|=|\lambda|,\ \mu\le\lambda} \F_\mu.
$$
In the classical case, $\F_{\leq \lambda}$ is spanned by certain products of
minors of the generic matrix. The {\it straightening law\/} of
Dubillet-Rota-Stein
[1974] shows that we get a basis if we choose only ``standard'' products of
these
types, and the successive quotients in the filtration
are the irreducible representations of ${\rm GL}(V)\times {\rm GL}(U)$.
The analogue
in our $\ZZ/2$-graded case is:

\proclaim Proposition 1.4. The subspaces $\F_{\le\lambda}$ define a
$\g$-invariant filtration on
$\S_{|\lambda|}(V\otimes U)$. The quotient $\F_{\le\lambda}/ \F_{<\lambda}$ is
isomorphic to $\bigwedge^\lambda V\otimes\bigwedge^\lambda
U=\S_{\lambda^\prime} V\otimes
\S_{\lambda^\prime} U$.\Box

There is also one element of each irreducible
representation which is easy to describe: the {\it highest weight
vector.} To speak of highest weight
vectors we must choose ordered bases $\{u_1,\dots,u_d\}$ and
$\{u'_1,\dots,u'_e\}$ of $U_0$ and $U_1$,
and ordered bases  $\{v_1,\dots,v_m\}$ and
$\{v'_1,\dots,v'_n\}$ of $V_0$ and $V_1$ respectively.

\proclaim Proposition 1.5.
Let
$\lambda = (\lambda_1 ,\ldots,\lambda_s)$
be a partition, and let
$w^1_i\in \bigwedge^{\lambda_i}(V)$ and
$w^2_i\in \bigwedge^{\lambda_i}(U)$ be the elements
$$
w_i^1 = \cases{ v_1\wedge\ldots\wedge v_{\lambda_i},&if
$\lambda_i\le m$;\cr
v_1\wedge\ldots\wedge v_m\wedge v_i'^{{(\lambda_i
-m)}},&otherwise.\cr}.
$$
$$
w_i^2 =\cases{ u_1\wedge\ldots\wedge u_{\lambda_i},&if
$\lambda_i\le d$;\cr
u_1\wedge\ldots\wedge u_d\wedge u_i^{\prime {(\lambda_i
-d)}},&otherwise.\cr}
$$
The element
$$
c_\lambda = \prod_{i=1}^s \rho_{\lambda_i} (w_i^1\otimes w_i^2 )
\quad \in\quad
\S(V\otimes U)
$$
is the  highest weight vector from the irreducible component
$\bigwedge^\lambda V\otimes \bigwedge^\lambda U
= \S_{\lambda'} V\otimes \S_{\lambda'} U$,
where
$\lambda'$
is the conjugate partition to $\lambda$.
\Box

We end this section with a desciption of a set of generators for
the representation $\bigwedge^\lambda V \otimes \bigwedge^\lambda U$.
We start with a double tableau, that is two sequences of tensors
$v_{i,1}\wedge\ldots\wedge v_{i,\lambda_i }in \bigwedge^{\lambda_i} V$ 
and 
$u_{i,1}\wedge\ldots\wedge u_{i,\lambda_i}\in\bigwedge^{\lambda_i}U$ 
($1\le i\le s$). We imagine that the elements
$v_{i,j}\in V$ correspond to the $i$-th row of the tableau $S$ of shape $\lambda$, 
and the elements $u_{i,j}\in U$ correspond to the $i$-th row of
another tableau $T$ of shape $\lambda$. We define
$$
\rho (S\otimes T)= 
\prod_{1\le i\le s} \rho_{\lambda_i}(v_{i,1}
\wedge\ldots\wedge 
v_{i,\lambda_i}\otimes
u_{i,1}\wedge\ldots\wedge u_{i,\lambda_i}).
$$

We think of $\lambda$ as a Ferrers diagram. If $S$ is a tableau of
shape $\lambda$ and $\sigma$ is a permutation of the boxes in $\lambda$
then $\sigma(S)$ is another tableau of shape $\lambda$ (here we write
$\sigma$ as a product of transpositions, and introduce a minus sign
whenever we interchange two elements of odd degree.)
Let $P(\lambda)$ be the group of permutations of the boxes in 
$\lambda$ that preserve the columns of $\lambda$.

\proclaim Proposition 1.6. The representation
$\bigwedge^\lambda V\otimes \bigwedge^\lambda U\subset R$
is generated by elements 
$$
\pi (S, T)=\sum_{\sigma\in P(\lambda )} \rho (\sigma S\otimes T ),
$$ 
or,
equivalently, by 
$$
\pi^\prime (S, T)=\sum_{\sigma\in P(\lambda )} \rho (S\otimes \sigma T)
$$
where $S$ and $T$ range over all tableaux of shape $\lambda$.

\noindent{\sl Proof.} We show that the $\pi(S,T)$ generate; the proof for $\pi'$ is
similar. 
Since $\rho(S,T)$ is antisymmetric in the elements appearing
in each row of $S$, the element 
$\sum_{\sigma\in P(\lambda)}\rho(\sigma(S), T)$ is 
the $\gl(V)$-linear projection of $\rho(S,T)$ to the 
$\bigwedge^\lambda V$-isotypic component of $R$.
By Corollary 1.2 we have 
$R=\oplus_\lambda \bigwedge^\lambda V\otimes \bigwedge^\lambda U$,
so this isotypic component is 
$\bigwedge^\lambda V\otimes \bigwedge^\lambda U\subset R$.\Box
\goodbreak
\beginsection 2. Proof of Theorem 1a.

In this section
$U$ and
$V$ are $\ZZ/2$-graded vector spaces of dimensions $(d,e)$
and $(m,n)$ respectively, and $\Phi$ is the generic map,
defined tautologically over $R=\S(V\otimes U)$.

We write $\Lambda(d,e)$ for the partition with
$e+1$ parts $((d+1)^e,d)=(d+1,\dots,d+1,d)$;
that is,
$\Lambda(d,e)$ corresponds to the Young diagram
that is a
$(d+1)\times (e+1)$ rectangle minus the box in the lower
right hand corner.
For example, $\Lambda(2,3)$ may be represented by the Young diagram
$$
\Lambda(2,3) =\Young (,,|,,|,,|,)\quad .
$$

For any
partition $\lambda$ we denote by $I_\lambda$ the
ideal in $R$ generated
by the representation $\bigwedge^\lambda V\otimes \bigwedge^\lambda U$.
With this notation, Theorem 1 takes the form:

\proclaim Theorem 1. The annihilator of the cokernel of $\Phi$
is equal to $I_{\Lambda(d,e)}$.

Theorem 1 implies that the representations
appearing in the annihilator of $\coker \Phi$ depend
only on the dimension of $U$, not the
dimension of $V$, as long as the dimension of $V_0$ is large
(Corollary 2.2), and we begin by proving this.
For the precise statement, we will use the following notation:
Let $V^\prime$ be another $\ZZ/2$-graded
vector space, and let
Let $\Phi^\prime$ be the generic map $R^\prime\otimes
V^\prime\rTo R^\prime\otimes U^*$
where $R^\prime = \S(V^\prime \otimes U )$.
If $V$ is a summand of $V'$, so that $V'=V\oplus W$,
then the ring $R=\S
(V\otimes U)$ can be identified with the subring of
$R^\prime$. We want to compare the
annihilators of the modules $\coker \Phi$ and
$\coker\Phi^\prime$.

\proclaim Proposition 2.1. If $V$ is a $\ZZ/2$-graded
 summand of $V'$
then
$\ann \coker\Phi=R\cap \ann \coker \Phi'$. More precisely,
\item{a)} $\coker\Phi$ is an $R$-submodule of $\coker \Phi^\prime $,
\item{b)} $\coker\Phi'$ is a quotient of
$(\coker \Phi)\otimes R'$.

\noindent{\sl Proof.} The first statement follows
easily from a) and b).

For the proof of a) and b) we may write
$V'=V\oplus W$ and we  make
use of the $\NN$-grading of $R'$ for which $V\otimes U$
has degree 0 and $W\otimes U$ has degree 1 (this grading
has nothing to do with the $\ZZ/2$-grading used elsewhere in
this paper!)
The map $\Phi^\prime$ is homogeneous of degree 0 if
we twist the summands of its source appropriately
$$
\Phi^\prime : R^\prime\otimes W(-1)\oplus R^\prime\otimes V
\rTo^{(\Phi'_1,\Phi'_0)}
R^\prime\otimes U^*,
$$
so we have an induced $\NN$ grading on $\coker\Phi^\prime$.
As $\Phi'_0=\Phi$, we see that $(\coker\Phi')_0=\coker \Phi$.
Since the elements of $R$ have degree 0, this is an $R$-submodule
as required for a).

For b) it suffices to note that $\coker \Phi'$ is obtained
from $(\coker\Phi)\otimes R'$ by factoring out the relations
corresponding to $W\otimes R'$.\Box

\proclaim Corollary 2.2. With notation $U,V,V',\Phi,\Phi'$
as above, suppose that $V$ is such that $I_\lambda\neq 0$
in $\S(V\otimes U)$. If $I_\lambda\subset \ann \coker \Phi$, then
$I_\lambda\subset \ann \coker\Phi^\prime$.

\noindent{\sl Proof.} 
The inclusion $R=\S(V\otimes U)\subset R'=\S(V^\prime \otimes U)$ carries
$\bigwedge^\lambda V\otimes \bigwedge^\lambda U$
into
$\bigwedge^\lambda V^\prime\otimes \bigwedge^\lambda U$.
The conclusion now follows from
Proposition 2.1 a) and b).
\Box
\medskip

\noindent{\sl Proof of Theorem 1.\/}
We first show that
$I_{\Lambda (d,e)}$ is contained in the annihilator of
$M:=\coker \Phi$. By Corollary 2.2, it is enough, given $U$, to
produce one nonzero element from
$
\bigwedge^{\Lambda(d,e)}V\otimes \bigwedge^{\Lambda(d,e)}
U\subset \S(V\otimes U)
$
that annihilates the cokernel of $\Phi$ for some space $V$.
By Corollary 2.2 it suffices to prove this result
in the case $m=d+1, n=0$, that is,  $dim\ V=(d+1,0)$.

Let $u_1 ,\ldots ,u_d$  be a basis of $U_0$, let
$u^\prime_1 ,\ldots ,u^\prime_e$ be
a basis of $U_1$,  and let $v_1 ,\ldots ,v_{d+1}$ be a basis
of $V=V_0$.
We denote the variables from the
$U_0\otimes V_0$ block by $x_{i,k}$ ($1\le i\le d,\ 1\le k\le d+1$),
and the variables from the $U_1\otimes V_0$ block by
$b_{j,k}$ ($1\le j\le e,\ 1\le k\le d+1$) thus:
$$\Phi =\pmatrix{
X\cr
B
}, X=(x_{i,k}), B=(b_{j,k}).
$$

Let
$$Z=Z_1 \cdot X(1,\ldots ,d\ |\ 1,\ldots ,d )$$
where $Z_1 =\prod_{j,k} b_{j,k}$
is the product of all the entries of $B$ and
$X(1,\ldots ,d\ |\ 1,\ldots ,d )$
is the $d\times d$ minor
of the matrix $X$ corresponding to the
first $d$ columns.

We now show that $Z$ annihilates $M$.
Indeed, by the classical
Fitting Lemma we know that $X(1,\ldots ,d\ |\ 1,\ldots ,d )$ annihilates
the even generic module $M/\Im\ (U_1^* )$. Thus
every basis element $u_k$ multiplied by
$X(1,\ldots ,d\ |\ 1,\ldots,d )$
can be expressed modulo the image of $\Phi$
as a linear combination of $u^\prime_1 ,\ldots ,u^\prime_e$ with
coefficients of positive degee in the
variables $b_{j,k}$. Since these variables are odd,
$Zu_k =0$ in $M$.

To see that $Z u_l^\prime$ is also 0 in $M$,
we use the classical Fitting Lemma  again on the
first $d$ columns of the matrix of $\Phi$ to see that for any
$1\leq i\leq d$ the element
$X(1,\ldots ,d\ |\ 1,\ldots ,d )u_i$
can be expressed,
modulo the image of $\Phi$, as
a linear combination of the $b_{s,t}u_j^\prime$.
On the other hand, if we multiply the last
column of the matrix of $\Phi$ by the product $Z_1'$ of all
the $b_{j,k}$ except for the $b_{l,d+1}$, we get
an expression for $Z_1 u_l^\prime$, modulo the
image of $\Phi$, as a linear
combination of $Z_1'u_1 \ldots ,Z_1'u_d$.
Thus $Zu'_l=X(1,\ldots ,d\ |\ 1,\ldots ,d )Z_1u'_l=0$ in $M$
as required.

Next we prove that the element $Z$ is a weight vector (not generally
a highest weight vector) and lies in
$\bigwedge^{\Lambda (d,e)} V\otimes \bigwedge^{\Lambda (d,e)} U$.
Indeed, the element $X(1,\ldots ,d\ |\ 1,\ldots ,d)$ is a weight vector in
$\bigwedge^d V\otimes\bigwedge^d U$.
The element $Z_1$ is a weight vector in the representation
$\bigwedge^{(d+1)^e}V\otimes \bigwedge^{(d+1)^e}U$.
The product is thus
contained in the ideal
$\F_{\le \Lambda(d,e)}$.
The element $Z$ has degree $(e+1)(d+1)-1$ but it involves only $d+1$
elements from $V=V_0$. By Proposition 1.5 its
weight can occur only in representations 
$\bigwedge^\lambda V\otimes \bigwedge^\lambda U\subset \S_\lambda(V\otimes U)$ 
with
$\lambda$ having all parts $\le d+1$. Since
$\Lambda (d,e)$ is the only partition $\lambda$ with
$\le e+1$ parts having
$|\lambda|=(d+1)(e+1)-1$
and each $\lambda_i\le d+1$, we are done.
This argument shows that
$I_{\Lambda(d,e)}$ is contained in the annihilator of the cokernel of
$\Phi$.

Now let $\mu$ be a partition not containing $\Lambda (d,e)$.
To complete the proof of Theorem 1, we must show that
the ideal $I_\mu$ does not
annihilate $M=\coker \Phi$ or, equivalently, that the highest
weight vector $c_\mu$ does not annihilate $M$.

Since $\mu$ does not contain
$\Lambda(d,e)$ it does not contain one of the extremal boxes of
$\Lambda(d,e)$.
By shifting the gradings of $V, U$ by $1$ we do not alter the
annihilator of the generic map, but we
change the notation so all partitions are changed to their
conjugates. Thus we may assume that $\mu_{e}\le d$.
By Corollary 2.2 we may further assume that $n=0$, so
that $V=V_0$, and that $m>>0$.
To prove the Theorem, we will do induction on $d$.

If $d=0$ we must show that the annihilator
of the cokernel of $\Phi$ is contained in
$I_{(1^e )}$; or equivalently that it contains
no $I_\lambda$ where $\lambda$ has fewer than $e$ parts.
Set
$$
Z_1=\prod_{1\le j\le e-1 ,1\le k\le m} b_{j,k};
$$
By Proposition 1.5,
$Z_1$ is the highest weight vector in
$\bigwedge^{m^{e-1}}V\otimes\bigwedge^{m^{e-1}}U$.
The element
$Z_1 u_e^\prime$ is not
in the image of $\Phi$
because the coefficient of
$u_e^\prime$ in any element from the image of $\Phi$ is in the
ideal generated by $b_{e,1},\ldots ,b_{e,m}$,
while $Z_1$ is not in this ideal.
Since $Z_1$ does not annihilate $M$, no $I_\lambda$
such that $\lambda$ has
$<e$ parts can annihilate $M$.

In case $d>0$ the matrix of $\Phi$ will contain an
even variable $x_{1,1}$. To complete the induction we will
invert this variable and use:

\proclaim Lemma 2.3.
\item{a)} Over the ring $R_1 =R[x_{1,1}^{-1}]$ the map $\Phi$ can be
reduced by row and column operations to the form
$$
\Phi^\prime\oplus\ id :\
(V^\prime\otimes_K R_1) \oplus R_1 \rightarrow
U^\prime {}^*\otimes_K \oplus R_1
$$
where $V$ is a $\ZZ/2$-graded vector space of dimension $(m-1,n)$ and
$U$ is a $\ZZ/2$-graded vector space
of dimension $(d-1,e)$. Moreover the ring $R^\prime$ generated over $K$ by
the entries of $\Phi^\prime$ is
isomorphic to $\S(V^\prime\otimes U^\prime )$ and $R_1$ is a flat
extension of $R^\prime$.
\item {b)} The localization of the ideal $I_\mu$ at $x_{1,1}$ is isomorphic
to the extension of the ideal
$J^\prime_\nu$ from $R^\prime$ where $\nu$ is the
partition obtained from $\mu$ by subtracting 1 from each nonzero part.

\noindent{\sl Proof of Lemma 2.3.}
Column and row reduction give the following formulas
for the entries of $\Phi^\prime$:
$$
x_{i,k}^\prime =x_{i,k}-{{x_{1,k} x_{i,1}\over x_{1,1}}},
\qquad
a_{i,l}^\prime
=a_{i,l}-{{a_{1,l} x_{i,1}\over
x_{1,1}}},
$$
$$
b_{j,k}^\prime =b_{j,k}-{{x_{1,k} b_{j,1}\over x_{1,1}}},
\qquad
y_{j,l}^\prime
=y_{j,l}-{{a_{1,l} b_{j,1}\over
x_{1,1}}}.
$$
Consequently
$$
R_1 = R^\prime [x_{1,1}, x_{1,1}^{-1}][x_{1,2},\ldots ,x_{1,m},\
a_{1,1},\ldots ,a_{1,n},\  x_{2,1},\ldots ,x_{d,1},\
b_{1,1},\ldots ,b_{e,1}]
$$
in the sense of $\ZZ/2$-graded algebras. This proves part a).

To prove part b) we first observe that the localization of the ideal $I_{(t
)}(\Phi )$ gives the ideal
$I_{(t-1)}(\Phi^\prime )$. Indeed, the ideal $I_{(t )}(\Phi )$
is generated by $\ZZ/2$-graded analogues of $t\times
t$ minors of $\Phi$. After localization it becomes the ideal $I_{(t
)}(\Phi^\prime \oplus id_{R_1} )$ generated by the $\ZZ/2$-graded analogues
of $t\times t$ minors of $\Phi^\prime\oplus
id_{R_1}$. Let us call the row and column of the matrix
 $\Phi^\prime\oplus id_{R_1}$ corresponding to the summand $R_1$ the
distinguished row and column respectively.
Every $\ZZ/2$-graded analogue of a $t\times t$ minor of
$\Phi^\prime\oplus id_{R_1}$ is either a $(t-1)\times (t-1)$minor of
$\Phi^\prime$ (in case it contains the distinguished row
and column), zero (if it contains the distinguished row but not the
distinguished column or vice versa), or a $t\times t$
minor of $\Phi^\prime$ if it does not contain the distinguished row nor
column.

To show that the result generalizes to an arbitrary partition
$\mu$ we order the bases so that the
distinguished row and column come first. We saw
in Proposition 6 that the highest weight vectors in $\bigwedge^\mu
V\otimes \bigwedge^\mu U$ are the products of minors of the matrix
$\Phi$ on
some initial subsets of rows and column of $\Phi$, so
after localization each factor will contain both the distinguished row and
the distinguished column of $\Phi^\prime\oplus
id_{R_1}$.
\Box

\bigskip
\noindent{\sl Completion of the Proof of Theorem 1.\/}
Now suppose that $d>0, n=0$. We may of course assume that $m\neq 0$,
so that the matrix of $\Phi$ contains the even variable $x_{1,1}$.
It is enough to prove that $I_\mu M\ne 0$ after inverting $x_{1,1}$.
The ideal
$I_{\mu}$ will localize to the ideal
$J^\prime_\nu$ where $\nu$ is equal to $\mu$ with all parts decreased by
$1$.  The graded vector space
$U$ of dimension $(d,e)$ will change to the $\ZZ/2$-graded vector space
$U^\prime$ of dimension
$(d-1,e)$. The desired conclusion follows by induction on $d$.\Box


\beginsection 3. Proof of Theorem 1b.

If $\phi: R^m\to R^d$ is a matrix representing a map of free modules
over a commutative ring then, as we noted in the introduction, there
are 
inclusions 
$ann(M)\cdot I_i(\phi)\subset I_{i+1}(\phi)$ for $0\leq i<d$ and thus, by induction,
$(ann(M)^d\subset I_d(\phi)$;
see for example Eisenbud[1995]. 
To
prove these inclusions one first notes that the cokernel of 
$\phi$ is the same as the cokernel of 
$$
\psi :V_0\otimes R\oplus U_0^*\otimes R\rightarrow U_0^*\otimes R
$$
where $\psi =(\phi, \, a\cdot Id )$. Thus $I_j (\phi )=I_j (\psi )$
and
$I_{j+1}(\phi )=I_{j+1}(\psi )\supset a\cdot I_j (\phi )$.
We will carry out the same
approach in $\ZZ/2$-graded case. 

In this section we work  with an arbitrary map
$$
\phi :V\otimes R\rightarrow U^*\otimes R
$$
of $\ZZ/2$-graded free modules over a $\ZZ/2$-graded commutative ring $R$.
The first step is to show that, just as in the classical case, 
the ideals $I_\lambda(\phi)$ depends only on the cokernel of $\phi$ and
on the number and degrees of the generators chosen.

\proclaim Lemma 3.1. 
If $\alpha :V'\otimes R\rightarrow V\otimes R$ then 
$$
I_\lambda (\phi\alpha)\subset I_\lambda (\phi).
$$
In particular, if
$\psi :V'\otimes R\rightarrow U^*\otimes R$ has the same cokernel
as $\phi$, then $I_\lambda (\phi)=I_\lambda (\psi)$.

\noindent{\sl Proof.} The second statement follows from the first because
each of the maps $\phi$ and $\psi$ factors through the other.

To prove the first statement, we use the notation of Proposition 1.6.
For any map $W\otimes R\to U^*\otimes R$, and any tableaux $S$ and $T$
of elements in $W$ and $U$, both of shape $\lambda$,
$\pi'_\psi(S,T)$ be the result of specializing the element $\pi'(S,T)$ defined
for the generic map $\Phi$ when $\Phi$ is specialized to $\psi$.
By Proposition 1.6 
it is enough to show that when $W=V'$ the element $\pi'_{\phi\alpha}(S,T)$
is in $I_\lambda(\phi)$. We have
$$
\rho_l(v_1'\wedge\dots v_l'\otimes u_1\wedge\dots u_l)
=
\sum_{i_1<\cdots<i_l}
\rho_l(v_1'\wedge\dots v_l'\otimes v_{i_1}^*\wedge\dots v_{i_l}^*)
\rho_l(v_{i_1}\wedge\dots v_{i_l}\otimes u_1\wedge\dots u_l)
$$
where $v_{1},\dots, v_{m+n}$ and 
$v_{1}^*,\dots, v_{m+n}^*$ are dual bases of $V$ and $V^*$.
Using this identity to rewrite the formula for
$\pi'_{\phi\alpha}(S',T)$, 
where $S'$ is a tableau of shape $\lambda$ with entries in $V'$
and $T$ is a tableau of shape $\lambda$ with entries in $U$, we see that
$\pi'_{\phi\alpha}(S',T)$ is a linear combination of elements of the form
$\pi'_\phi(S,T)$, where 
$S$ is a tableau of shape $\lambda$ with entries in $V$.
\Box

Lemma 3.1 implies in particular that two presentations of the same module
with the same numbers of even and odd generators have the same 
ideals $I_\lambda(\phi)$. Similar arguments show that we can allow
for presentations with different numbers of generators as long 
as we change the partitions suitably: if we add $d'$ even and $e'$ odd 
generators, then we have to expand $\lambda$ by adding $d'$ columns
of length equal to the length of the first column  and 
$e'$ rows of length equal to the length of the first row of the
resulting partition (or vice versa). In this sense 
the ideals $I_\lambda(\phi)$ depend only on the cokernel of 
$\phi$.

The main result of this section is the following

\proclaim Theorem 3.2. Let $R, U, V, \phi$ be as in 
the beginning of the introduction, and let $M=\coker \phi$.
\item {a)} Let $s$ be an integer, $0\le s\le d-1$. 
If $x_1 ,\ldots ,x_{e+1}\in Ann_R M$, 
then 
$$
x_1\ldots x_{e+1}I_{\Lambda (s,e)}(\phi )\subset I_{\Lambda (s+1,e)}(\phi ).
$$
\item{b)} If $x_1 ,\ldots ,x_e \in Ann_R M$, then
$x_1\ldots x_e\in I_{\Lambda (0,e)}(\phi)$.

As in the classical case we derive:

\proclaim Corollary 3.3. Let $M$ be a 
$\ZZ/2$-graded module over a $\ZZ/2$-graded ring $R$,
with the presentation
$\phi :V\otimes R\rightarrow U^*\otimes R$. 
Assume that $dim\ U = (d,e), dim\ V=(m,n)$. Let
$x_1 ,\ldots,x_{(d+1)(e+1)-1}$ be  homogeneous elements from
$Ann_R M$. Then $x_1 \ldots x_{(d+1)(e+1)-1}\in I_{\Lambda (d,e)}(\phi )$.
\Box

\noindent {\sl Proof of Theorem 3.2.} 
We begin with part b).  We work with a presentation 
$(\phi ,\psi ):V\otimes R\oplus W\otimes R\rightarrow U^*\otimes R$ 
where $W$ is a
$\ZZ/2$-graded vector space of dimension $e$ with the the $i$-th
generator $w_i$ going to $x_i$ times the $i$-th generator $u_i$ of
$U^*$. The parity of the generators of $W$ is adjusted so $\psi$ is of
degree $0$. Now taking a double tableau $ (S, T)$ of the shape $(1^e)$
with $w_i$ and $u_i$ in the i-th row, and applying the definition
above, we see that the generator $\pi (S, T)$ is just $x_1\ldots x_e$.

To prove part a) we distinguish two cases. In the case $s<d-1$ we use
the presentation $(\phi ,\psi ) :V\otimes R\oplus W\otimes
R\rightarrow U^*\otimes R$ where $W$ is a $\ZZ/2$-graded vector space
of dimension $e+1$ with the the $i$-th generator $w_i$ going to $x_i$
times the $s+i+2$-nd generator $u_i$ of $U^*$.  The parity of the
generators of $W$ is adjusted so $\psi$ is of degree $0$. We can
assume without a loss of generality that $\phi :=\Phi$ is
generic. Then it is enough to prove that $c_{\Lambda (s,e)}x_1 \ldots
x_{e+1}\in I_{\Lambda (s+1,e)}$, where $c_{\Lambda (s,e)}$ is the
highest weight vector defined as in Proposition 1.5.

We pick a tableau $(S, T)$ of shape $\Lambda (s+1,e )$ as follows. The
entries $v_{i,j}$, $u_{i,j}$ in the $i$-th row are the same as in the
canonical tableau, excepr the last ones. The last entry in the tableau
$v$ in the $i$-th row is $w_i$, and the last entry in the $i$-th row
is $u_{s+2}$. The element $\pi^\prime (S, T)$ is easily seen to be
$c_{\Lambda (s,e)}x_1 \ldots x_{e+1}$.

In the case $s=d-1$ we use the presentation $(\phi ,\psi ) :V\otimes
R\oplus W\otimes R\rightarrow U^*\otimes R$ where $W$ is a
$\ZZ/2$-graded vector space of dimension $e+1$ with the the $i$-th
generator $w_i$ going to $x_i$ times the $d+i$-th generator $u_i$ of
$U^*$ for $1\le i\le e$ and $w_{e+1}$ goes to $x_{e+1}$ times $u_d$.
The parity of the generators of $W$ is adjusted so $\psi$ is of degree
$0$. We can assume without a loss of generality that $\phi :=\Phi$ is
generic. Then it is enough to prove that $c_{\Lambda (d-1,e)}x_1
\ldots x_{e+1}\in I_{\Lambda (d,e)}$, where $c_{\Lambda (d-1,e)}$ is
the canonical tableau.

We pick a tableau $(S, T)$ of shape $\Lambda (s+1,e )$ as follows. The
entries $v_{i,j}$, $u_{i,j}$ in the $i$-th row are the same as in the
canonical tableau, except the last ones. The last entry in the tableau
$v$ in the $i$-th row is $w_i$, and the last entry in the $i$-th row
is $u_{d+i}$ for $1\le i\le e$, and $u_d$ for $e+1$-st row. The
element $\pi^\prime (S, T)$ is easily seen to be $c_{\Lambda
(d-1,e)}x_1 \ldots x_{e+1}$.$\Box$
\goodbreak

\beginsection 4. The resolution of generic $\ZZ/2$-graded module.

In this section we work over the generic ring $R=\S$ as in the
introduction, and we conjecture the form of a
minimal free resolution over $R$ of the cokernel $C$ of the
generic map $\Phi$. This
resolution is a natural generalization of the one constructed in [B-E] in
the commutative case. We work over a field $K$
of characteristic $0$. We define some
$\ZZ/2$-graded free $R$-modules $\FF_i$ as follows:
$$
{\FF}_0 = U^*\otimes R,\qquad {\FF}_1 = V\otimes R
$$
$$
{\FF}_i =\oplus_{|\alpha |+|\beta | = i-2}  \S_{\Theta (d,e,\alpha
,\beta )}V\otimes \S_{\Lambda (d,e,\alpha
,\beta )}U\otimes R
$$
where $\Lambda (d,e,\alpha ,\beta ) = (d+1+\beta_1 ,d+1+\beta_2,\ldots
d+1+\beta_e ,e,\alpha^\prime_1 ,\ldots
,\alpha^\prime_s )$, $\Theta (d,e,\alpha ,\beta )= (d+1+\alpha_1
,d+1+\alpha_2,\ldots d+1+\alpha_e
,e+1,\beta^\prime_1 ,\ldots ,\beta^\prime_s )$, and we sum over all pairs
of partitions $\alpha ,\beta$ with at
most $e$ parts.

\proclaim Conjecture 4.1. There exists an
 equivariant differential $d_i :{\FF}_i\rightarrow {\FF}_{i-1}$,
linear for
$i\ge 3, i=1$ and of degree
$|\Lambda (d,e)|$ for $i=2$, that makes
${\FF}_\bullet$ into a minimal
$R$-free $\ZZ/2$-graded resolution of $C$.

In the even case ($U_1=V_1=0$) the desired complex is the
Buchsbaum-Rim resolution (see for example Buchsbaum and Eisenbud [1975]).
We have
checked the conjecture computationally, using Macaulay2,
in a few more cases.

\beginsection Appendix. Comments on the action of $\g$

It may at first be surprising that the generators given in 
Examples 1-3 of the Introduction are permuted by the action of $\g$,
so we pause to make the action explicit in Example 3 (the
other cases are similar and simpler). 

When we think of $R=\S(V\otimes U)$ as a $\g$ module, we think of
$\g$ acting on the left. But we may identify $U\otimes V$
with ${\rm Hom}(V,U^*)^*={\rm Hom}(U^*,V)$, and thus identify $R$ with the
coordinate ring of the space ${\rm Hom}(V,U^*)$. In this identification
it is natural to think of the Lie algebra $\g=\gl(V)\times\gl(U)$
as $\gl(V)\times\gl(U^*)$, with the $\gl(U^*)$ acting on the right.
To make this identification, we use the {\it supertranspose\/} which
is the anti-isomorphism 
$$
\gl(U)\to \gl(U^*);\quad 
\left(\matrix{U_{0,0}&U_{0,1}\cr
U_{1,0}&U_{1,1}}\right)
\mapsto
\left(\matrix{U_{0,0}^t&U_{1,0}^t\cr
-U_{0,1}^t&U_{1,1}^t}\right)
$$

Now consider the case presented in Example 3 of the introduction, 
whose notation we use.
To see the action of $\g$ let us act by two elements of the Lie algebra on
the element $axy$. First we act with
the element $v_{0,1}$ from $\gl(V)$ changing an odd element to an even one.
We get a sum of terms each of
which is $axy$ with one changed factor; we can replace $a$ by $x$ or $y$ by
$b$. Thus we get terms $xxy$ and
$axb$. The first comes with positive sign ($v_{0,1}$ acts from the left, we
replaced first factor, so there
are no switches), the second term comes with negative sign (we replaced $y$
by $b$, so had to switch $v_{0,1}$
with $a$). Thus we get $x(xy-ab)$.

Let us also act on $axy$ by the Lie algebra element $u_{1,0}$
from $\gl(U^*)$ exchanging an even element with an odd one.
We get a sum of terms where each term is $axy$ with one factor
changed; we can
change $x$ to $b$ and $a$ to $y$. We get terms $yxy$ and $aby$. Both come
without sign as $u_{1,0}$ acts from
the right and $x,y$ have even degree. Thus we get $(xy+ab)y$.

\beginsection References.

{\frenchspacing


\noindent A. Berele and A. Regev: Hook Young diagrams with applications
to combinatorics and to representations of Lie superalgebras. Adv. in Math. 64
(1987) 118--175.
\smallskip

\noindent D. A. Buchsbaum and D. Eisenbud: Generic free
resolutions and a family of generically perfect
ideals. Advances in Math. 18  (1975), no. 3, 245-301.
\smallskip

\noindent C. DeConcini, D. Eisenbud and C. Proces: Young
diagrams and determinantal varieties. Invent.
Math. 56  (1980) 129--165.
\smallskip

\noindent P. Doubilet, G.-C. Rota and Joel Stein:
On the foundations of combinatorial theory. IX.
Combinatorial methods in invariant theory.
Studies in Appl. Math. 53 (1974), 185--216.
\smallskip

\noindent D. Eisenbud: {\it Commutative Algebra with a View Toward
Algebraic Geometry}. Springer-Verlag, New York, 1995.
\smallskip

\noindent D. Eisenbud, S. Popescu and S. Yuzvinsky: Hyperplane arrangement
cohomology and monomials in the exterior algebra. Preprint, 2000.
\smallskip

\noindent D. Eisenbud and F.-O. Schreyer: Free resolutions and
sheaf cohomology over exterior algebras. Preprint, 2000.
\smallskip

\noindent H.~Fitting: Die Determinantenideale eines Moduls. Jahresbericht
der Deutschen Math.-Vereinigung 46 (1936) 195-229.
\smallskip

\noindent D.~Grayson, M.~Stillman: {\it Macaulay2}, a
software system devoted to supporting research in algebraic geometry
and commutative algebra.  Contact the authors, or download from\hfill\break
{\tt http://www.math.uiuc.edu/Macaulay2} .
\smallskip

\noindent M. Green: The Eisenbud-Koh-Stillman conjecture on linear syzygies.
Invent. Math. 136 (1999)
411--418.
\smallskip

\noindent I. G. Macdonald: {\it Symmetric functions and Hall polynomials.}
Second edition. With contributions by A.
Zelevinsky. Oxford Mathematical Monographs,
Oxford University Press, New York, 1995.
\smallskip

}
\bigskip\bigskip
\noindent {\bf Author Addresses:}
\smallskip
\noindent{David Eisenbud}\par
\noindent{Department of Mathematics, University of California, Berkeley,
Berkeley CA 94720}\par
\noindent{de@msri.org}
\smallskip
\noindent{Jerzy Weyman}\par
\noindent Department of Mathematics, Northeastern  University,
Boston MA 02115\par
\noindent j.weyman@neu.edu

\bye